\documentclass[12pt]{amsart}
\usepackage{amsmath,amssymb,latexsym}

\def\Ind#1#2{#1\setbox0=\hbox{$#1x$}\kern\wd0\hbox to 0pt{\hss$#1\mid$\hss}
\lower.9\ht0\hbox to 0pt{\hss$#1\smile$\hss}\kern\wd0}

\def\Notind#1#2{#1\setbox0=\hbox{$#1x$}\kern\wd0\hbox to 0pt{\mathchardef
\nn="3236\hss$#1\nn$\kern1.4\wd0\hss}\hbox to 0pt{\hss$#1\mid$\hss}\lower.9\ht0
\hbox to 0pt{\hss$#1\smile$\hss}\kern\wd0}

\def\Var{\mathrm{Var}}
\def\var{\mathrm{var}}
\def\ac{\mathrm{ac}}
\def\Mc{$\mathfrak M$c}

\newtheorem{theorem}{Theorem}[section]
\newtheorem{fact}[theorem]{Fact}
\newtheorem{cor}[theorem]{Corollary}
\newtheorem{conj}[theorem]{Conjecture}
\newtheorem{lmm}[theorem]{Lemma}
\newtheorem*{claim}{Claim}
\theoremstyle{definition}
\newtheorem{remark}[theorem]{Remark}
\newtheorem{definition}[theorem]{Definition}
\newtheorem{expl}[theorem]{Example}
\newtheorem{quest}[theorem]{Problem}
\newtheorem*{nota}{Notation}
\def\bsp{\begin{expl}}
\def\ebsp{\end{expl}}
\def\frag{\begin{quest}}
\def\efrag{\end{quest}}
\def\beh{\begin{claim}}
\def\ebeh{\end{claim}}
\def\defn{\begin{definition}\upshape}
\def\edefn{\end{definition}}
\def\satz{\begin{theorem}}
\def\esatz{\end{theorem}}
\def\tats{\begin{fact}}
\def\etats{\end{fact}}
\def\kor{\begin{cor}}
\def\ekor{\end{cor}}
\def\lem{\begin{lmm}}
\def\elem{\end{lmm}}
\def\bem{\begin{remark}}
\def\ebem{\end{remark}}
\def\verm{\begin{conj}}
\def\everm{\end{conj}}
\def\bew{\begin{proof}}
\def\bewbeh{\begin{proof}[Proof of Claim}
\def\ebew{\end{proof}}

\begin{document}
\title{Largeness and equational probability in groups}
\author{Khaled Jaber}
\address{Khaled Jaber, Department of Mathematics, Faculty of Sciences,
Lebanese University, Lebanon}
\email{kjaber@ul.edu.lb}
\author{Frank O. Wagner}
\address{Frank O.~Wagner, Universit\'e de Lyon; Universit\'e
Lyon 1; CNRS UMR 5208, Institut Camille Jordan, 21 avenue Claude Bernard,
69622 Villeurbanne-cedex, France}
\thanks{The second author was partially supported by the ANR-DFG project AAPG2019 GeoMod.}
\email{wagner@math.univ-lyon1.fr}

\date{\today}
\begin{abstract}We define $k$-genericity and $k$-largeness for a subset of a group, and determine the value of $k$ for which a $k$-large subset of $G^n$ is already the whole of $G^n$, for various equationally defined subsets. We link this with the inner measure of the set of solutions of an equation in a group, leading to new results and/or proofs in equational probabilistic group theory.\end{abstract} 
\maketitle

\section{Introduction}
In probabilistic group theory we are interested in what proportion of (tuples of) elements of a group have a particular property; if this property is given by an equation, we talk about {\em equational probability}. In \cite{JW00} a notion of {\em largeness} was introduced for a subset of a group, and it was shown that certain equational properties of a group hold everywhere as soon as they hold largely. In this paper, we shall introduce a quantitative version of largeness, and deduce some results in equational probabilistic group theory.

Throughout this paper, $G$ will be a group and $\mu$ a left-invariant probability measure on some algebra of subsets of $G$.
\bsp\begin{enumerate}\item $G$ finite, $\mu$ the counting measure.
\item $G_1$ a group, $\mu_1$ a left-invariant measure on $G_1$, and $G=G_1^n$ with the product measure $\mu=\mu_1^n$.
\item More generally, $G_1$ a group, $G\le G_1^n$ and $\mu$ a left-invariant measure on $G$.
\item $G$ arbitrary and the measure algebra reduced to $\{\emptyset,G\}$. While this set-up trivialises the probability statements, the largeness results remain meaningful.\end{enumerate}\ebsp
If $X$ is a measurable subset of $G$ we can interpret $\mu(X)$ as the probability that a random element of $G$ lies in $X$. If $H$ is another group, $f:G\to H$ is a function and $c\in H$ some constant, we put $\mu(f(x)=c)=\mu(\{g\in G:f(g)=c\})$.
\bsp\label{b:f} Let $G_1$ be a group, $G\le G_1^n$ a subgroup, $\bar g\in G_1^m$ constants, and $w(\bar x,\bar y)$ a word in $\bar x\bar y $ and their inverses, with $|\bar x|=n$ and $|\bar y|=m$. Then $w(\bar x,\bar g)$ induces a function from $G$ to $G_1$.\ebsp

We shall now list some known results, starting with Frobenius in 1895.
\begin{fact}\label{f:Fact} Let $G$ be a finite group.\begin{itemize}
\item {\bf Frobenius 1895 \cite{Fr95}} If $n$ divides $|G|$ then the number of solutions of $x^n=1$ is a multiple of $n$. In particular, $\mu(x^n=1)\ge\frac{n}{|G|}$.
\item{\bf Miller 1907 \cite{Mi07}} If $G$ is non-abelian, then $\mu(x^2=1)\le\frac34$.
\item{\bf Laffey 1976 \cite{La76}} If $G$ is a $3$-group not of exponent $3$ then $\mu(x^3=1)\le\frac79$.
\item{\bf Laffey 1976 \cite{La76a}} If $p$ is prime and divides $|G|$, but $G$ is not a $p$-group, then $\mu(x^p=1)\le\frac p{p+1}$.
\item{\bf Laffey 1979 \cite{La79}} If $G$ is not a $2$-group, then $\mu(x^4=1)\le\frac89$.
\item{\bf Iiyori, Yamaki 1991 \cite{IY91}} If $n$ divides $|G|$ and $X=\{g\in G:g^n=1\}$ has cardinality $n$, then $X$ forms a subgroup of $G$.
\item {\bf Erd\H os, Turan, 1968 \cite{ET68}} If $k(G)$ is the number of conjugacy classes in $G$, then $\mu([x,y]=1)=\frac{k(G)}{|G|}$.
\item {\bf Joseph 1977 \cite{Jo77}, Gustafson 1973 \cite{Gu73}} If $G$ is non-abelian, then $\mu([x,y]=1)\le\frac58$.
\item {\bf Neumann, 1989 \cite{Ne89}} For any real $r>0$ there are  $n_1(r)$ and $n_2(r)$ such that if $\mu([x,y]=1)\ge r$ then $G$ contains normal subgroups $H\le K$ such that $K/H$ is abelian, $|G:K|\le n_1(r)$ and $|H|\le n_2(r)$.
\item {\bf Barry, MacHale, N\'\i\ Sh\'e, 2006 \cite{BMS06}} If $\mu([x,y]=1)>\frac13$ then $G$ is supersoluble.
\item {\bf  Heffernan, MacHale, N\'\i\ Sh\'e, 2014 \cite{HMS14}} If $\mu([x,y]=1)>\frac7{24}$ then $G$ is
metabelian. If $\mu([x,y]=1)>\frac{83}{675}$ then $G$ is abelian-by-nilpotent.\end{itemize}
\end{fact}

In Section \ref{s:Large} we shall introduce largeness and prove the main connection between largeness and measure, Lemma \ref{l:mu}, which will be used throughout the rest of the paper. Section \ref{s:FC} deals with central elements, or more generally FC and BFC groups.
We shall treat equations of the form $x^n=c$ for arbitrary $c$ in Section \ref{s:Burnside}, recovering Miller's result for $n=2$, and a weaker bound than Laffey for $n=3$ (namely $\frac67$). In Section \ref{s:[]} we shall consider commutator equations; while our methods allow us to deal with more complicated commutators, they are too general to obtain the bounds from Fact \ref{f:Fact}.
Section \ref{s:Nilp} deals with nilpotent groups via linearisation, and the short Section \ref{s:Auto} places Sherman's autocommutativity degree in our context.

\begin{nota} We shall write $x^y=y^{-1}xy$, $x^{-y}=(x^{-1})^y=y^{-1}x^{-1}y$ and $[x,y]=x^{-1}y^{-1}xy=y^{-x}y=x^{-1}x^y$.
\end{nota}

\section{Largeness and Probability}\label{s:Large}
The following notion of largeness was introduced in \cite{JW00}.
\defn If $X\subseteq G$, we say that $X$ is {\em $k$-large} in $G$ if the intersection of any $k$ left translates of $X$ is non-empty, and $X$ is {\em $k$-generic} in $G$ if $k$ left translates of $X$ cover $G$.
%The {\em genericity index} $\ig(X)$ of $X$ is the minimal $k$ such that $X$ is 
%$k$-generic; the {\em largeness index} $\il(X)$ of $X$ is the maximal $k$ such that $X$ is $k$-large.  
A subset $X$ is {\em large} if it is $k$-large for all $k$; it is {\em generic} if it is $k$-generic for some $k$.\edefn
Of course, analogous notions exist for right and two-sided genericity/largeness. Both genericity and largeness are notions of prominence, increasing with $k$ for largeness and decreasing with $k$ for genericity. Clearly, if $X\subseteq G$ and $X$ is ($k$-)large/generic, so is any left or right translate or superset of $X$. Largeness and genericity are co-complementary:
\lem Let $X\subseteq G$. Then $X$ is $1$-large if and only if $X\not=\emptyset$, and $X$ is $1$-generic if and only if $X=G$. More generally, $X$ is $k$-large if and only if $G\setminus X$ is not $k$-generic.  Finally, $X$ is $k$-generic/large if and only if $X\cap Y\not=\emptyset$ for all $k$-large/generic $Y\subseteq G$.\elem
\bew We only show the last assertion. If $X$ is not $k$-generic/large, then $Y:=G\setminus X$ is $k$-large/generic, and $X\cap Y=\emptyset$. Conversely, if $X$ is $k$-generic, say $G=\bigcup_{i<k}g_iX$, and $Y$ is $k$-large, then
$$\begin{aligned}\emptyset&\not=\bigcap_{i<k}g_iY=G\cap\bigcap_{i<k}g_iY
=\bigcup_{i<k}g_iX\cap\bigcap_{i<k}g_iY\\
&=\bigcup_{i<k}(g_iX\cap\bigcap_{i<k}g_iY)
\subseteq\bigcup_{i<k}(g_iX\cap g_iY)=\bigcup_{i<k}g_i(X\cap Y).\end{aligned}$$
Thus $X\cap Y\not=\emptyset$.
\ebew
\bem If $\phi:G\to H$ is an epimorphism and $X\subseteq G$ is ($k$-)large/generic, so is $\phi(X)\subseteq H$.
Conversely, if $Y\subseteq H$ is ($k$-)large/generic in $H$, so is $\phi^{-1}[X]$ in $G$.

In particular, if $X\subseteq G\times H$ is ($k$-)large/generic, so are the projections to each coordinate. Conversely, if $X\subseteq G$ and $Y\subseteq H$ are ($k$-)large, so is $X\times Y\subseteq G\times H$; if $X$ is $k$-generic and $Y$ is $\ell$-generic, $X\times Y$ is $k\ell$-generic.\ebem
\lem Suppose $X$ is $k\ell$-large in $G$ and $H\le G$ is a subgroup of index $k$. Then $X\cap H$ is $\ell$-large in $H$.\elem
\bew Let $(g_i:i<k)$ be coset representatives of $H$ in $G$, and consider $(h_j:j<\ell)$ in $H$. By $k\ell$-largeness of $X$ in $G$ there is $x\in\bigcap_{i<k,\ j<\ell}g_ih_jX$. As $\bigcup_{i<k}g_iH=G$, there is $i_0<k$ with $x\in g_{i_0}H$. But then
$$g_{i_0}^{-1}x\in H\cap\bigcap_{i<k,\ j<\ell}g_{i_0}^{-1}g_ih_jX\subseteq H\cap\bigcap_{j<\ell}h_jX
=\bigcap_{j<\ell}h_j(X\cap H),$$
so $X\cap H$ is $\ell$-large.\ebew

The link between largeness and probability is given by the following lemma, which will be used throughout the paper. Recall that the {\em inner measure} of an arbitrary subset $X$ of a measurable group $G$ is
$$\mu_*(X)=\sup\{\mu(Y):Y\subseteq X\mbox{ measurable}\},$$
and the {\em outer measure} is given by
$$\mu^*(X)=\inf\{\mu(Y):Y\supseteq X\mbox{ measurable}\}.$$
Clearly the inner measure is superadditive, the outer measure is subadditive, and $\mu_*(X)+\mu^*(G\setminus X)=1$.
\lem\label{l:mu} If $X$ is $k$-generic in $G$, then $\mu^*(X)\ge\frac1k$. If $\mu_*(X)>1-\frac1k$ then $X$ is $k$-large in $G$.\elem
\bew If $X$ is $k$-generic there are $g_1,\ldots,g_k$ in $G$ with $G=\bigcup_{i\le k}g_iX$. Hence
$$1=\mu^*(G)=\mu^*(\bigcup_{i\le k}g_iX)\le\sum_{i\le k}\mu^*(g_iX)=k\,\mu^*(X)$$
by left invariance, whence $\mu^*(X)\ge\frac1k$.

Now if $X$ is not $k$-large, its complement is $k$-generic, so $\mu^*(G\setminus X)\ge\frac1k$. But then $\mu_*(X)\le1-\frac1k$.\ebew
These bounds are strict, as we can take $X$ a subgroup of index $k$ (resp.\ its complement).
\bem For any group $G$ the set $(G\times\{1\})\cup(\{1\}\times G)$ is $2$-large in $G^2$; if $G$ is infinite it is of measure $0$.\ebem

We shall now prove some results about finite groups, which owing to their non-linearity do not generalise easily to the measurable context.
\bem Let $G$ be a finite group of order $n$, and $X\subseteq G$ a non-empty proper subset of size $m$. Then $X$ is $(n-m+1)$-generic and at most $m$-large, since we can form the union of $X$ with $n-m$ translates of $X$ to cover all the $n-m$ points of $G\setminus X$, and we can intersect $X$ with $m$ translates of $X$ to remove all $m$ points of $X$.\ebem
\satz Let $G$ be a finite group of order $n$, and $X\subseteq G$ a non-empty proper subset of size $m$. If $m>n-\frac12-\sqrt{n-\frac34}$, then $X$ is $2$-generic. Hence if $m<\frac12+\sqrt{n-\frac34}$ then $X$ is not $2$-large.\esatz
\bew If $m>n-\frac12-\sqrt{n-\frac34}$, then
$$\textstyle n-\frac34>(n-m-\frac12)^2=(n-m)(n-m-1)+\frac14.$$
Put $Z=\{xy^{-1}:x,y\in G\setminus X\}$. Then
$$|Z|\le(n-m)(n-m-1)+1<n,$$
so there is $g\in G\setminus Z$. But if $h\in G\setminus(X\cup gX)$, then $h,g^{-1}h\in G\setminus X$, and $g=h(g^{-1}h)^{-1}\in Z$, a contradiction. Thus $G=X\cup gX$ and $X$ is $2$-generic.

The second assertion follows by taking complements.\ebew
\satz Let $G$ be a finite group of order $n$. If the exponent of $G$ does not divide $\ell$ then $\mu(x^\ell=1)\le1-\frac1{\sqrt{2n}}$, where $\mu$ is the counting measure.\esatz
\bew Put $X=\{g\in G:g^\ell=1\}$, of size $m<n$, and take any $g\in G\setminus X$. Note that $X\cap gX\cap C_G(g)$ is empty, as otherwise there would be $y\in C_G(g)$ with $y^\ell=1=(gy)^\ell$, whence $g^\ell=1$ and $g\in X$.

Thus $|C_G(g)|\le 2\,|G\setminus X|$. Moreover $g^G\cap X=\emptyset$, and 
$$|G|/|C_G(g)|=|g^G|\le|G\setminus X|.$$ 
Thus $n=|G|\le 2\,|G\setminus X|^2$ and $\sqrt{\frac n2}\le n-m$, whence
$${\qquad\qquad\qquad}\mu(x^\ell=1)=\frac mn\le \frac{n-\sqrt{\frac n2}}n=1-\frac1{\sqrt{2n}}.\qedhere$$
\ebew

 \defn Let $f:G\to H$ be a function, and $c\in H$. The equation $f(x)=c$ is {\em $k$-largely satisfied} in $G$ if $\{g\in G:f(g)=c\}$ is $k$-large in $G$. By abuse of notation, if $G=G_1^n$ and $x=(x_1,\ldots,x_n)$, we shall also say that $f(x_1,\ldots,x_n)=c$ is $k$-largely satisfied in $G_1$.\edefn

\section{FC-Groups}\label{s:FC}
In this section we shall work in the set-up of Example \ref{b:f}: $G_1$ will be a group, $G\le G_1^n$, $w(\bar x,\bar y)$ a word in $\bar x\bar y$ and their inverses with $n=|\bar x|$ and $m=|\bar y|$, $\bar g\in G_1^m$ and $c\in G_1$ constants, and $f(\bar x)=w(\bar x,\bar g)$.\medskip

Recall that a group is $FC$ if the centraliser of any element has finite index; it is $BFC$ if the index is bounded independently of the element.

We shall first need a preparatory lemma.
For two tuples $\bar g=(g_i:i<k)$ and $\bar g'=(g'_i:i<k)$ in $G_1^k$ we shall put $\bar g^{-1}=(g_i^{-1}:i<k)$ and $\bar g\cdot\bar g'=(g_ig'_i:i<k)$.

\lem\label{l:3} Suppose $\bar g,\bar g'\in G_1^m$ and $\bar h,\bar h'\in G_1^n$ are such that all elements from $\bar g\bar h$ commute with all elements from $\bar g'\bar h'$. Then 
$$w(\bar h\cdot\bar h',\bar g\cdot\bar g')=w(\bar h,\bar g)\,w(\bar h',\bar g').$$\elem
\bew Obvious.\ebew

\satz Let $G_1$ be an $FC$-group. If the equation $w(\bar x,\bar g)=c$ is largely satisfied in $G$ then it is identically satisfied in $G$.\esatz
\bew Consider $\bar h\in G$, and $C=C_{G_1}(\bar g,\bar h)$, a subgroup of finite index in $G_1$. Put $H=C^n\cap G$, a subgroup of finite index in $G$, and $X=\{\bar h'\in G:w(\bar h',\bar g)=c\}$. Then $X\cap \bar h^{-1}X\cap H$ is large in $H$, whence non-empty. So there is $\bar x\in H$ with
$$w(\bar 1,\bar g)\,w(\bar x,\bar 1)=w(\bar x,\bar g)=c=w(\bar h\cdot\bar x,\bar g)=w(\bar h,\bar g)\,w(\bar x,\bar 1).$$
Hence $w(\bar h,\bar g)=w(\bar 1,\bar g)$ for all $\bar h\in G$, and $w(\bar 1,\bar g)=w(\bar x,\bar g)=c$.\ebew

For a $BFC$-group, we can bound the degree of largeness needed:
\satz\label{BFC} Suppose every centraliser of a single element has index at most $k$ in $G_1$. If the equation $w(\bar x,\bar g)=c$ is $2k^{n^2+mn}$-largely satisfied in $G$ then it is identically satisfied in $G$.\esatz
\bew In the notation of the previous proof, $C=C_{G_1}(\bar g,\bar h)$ has index at most $k^{n+m}$ in $G_1$, so 
$$|G:H|=|G:G\cap C^n|\le|G_1^n:C^n|=|G_1:C|^n\le (k^{n+m})^n=k^{n^2+mn}.$$
Now $2k^{n^2+mn}$-largeness of $X$ in $G$ implies $k^{n^2+mn}$-largeness of $X\cap \bar h^{-1}X$ in $G$, whence $1$-largeness of $X\cap \bar h^{-1}X\cap H$ in $H$. So we can find the $\bar x$ required to finish the proof.\ebew
\kor\label{c:6} Suppose every centraliser of a single element has index at most $k$ in $G_1$. If $w(\bar x,\bar g)=c$ is not an identity on $G$, then 
$$\mu_*(w(\bar x,\bar g)=c)\le1-\frac1{2k^{n^2+mn}}.$$\ekor
\bew If $\mu_*(w(\bar x,\bar g)=c)>1-\frac1{2k^{n^2+mn}}$, then $\{\bar x\in G:w(\bar x,\bar g)=c\}$ is $2k^{n^2+mn}$-large in $G$ by Lemma \ref{l:mu}, and $w(\bar x,\bar g)=c$ is identically satisfied in $G$ by Theorem \ref{BFC}.\ebew
\bem\label{b:exp} This holds in particular for the equation $x^\ell=c$, with $n=1$ and $m=0$.\ebem

If the group is central-by-finite, the largeness needed does not depend on the number of parameters.
\kor\label{k:zentrum} Suppose $Z(G_1)$ has index $k$ in $G_1$. If the equation $w(\bar x,\bar g)=c$ is $2k^n$-largely satisfied in $G$ then it is identically satisfied in $G$.\ekor
\bew $H=G\cap Z(G_1)^n$ has index at most $k^n$ in $G$. We finish as in Theorem \ref{BFC}.\ebew
\kor If $|G_1:Z(G_1)|\le k$ and $w(\bar x,\bar g)=c$ is not an identity in $G$, then $\mu_*(w(\bar x,\bar g)=1)\le1-\frac1{2 k^n}$.\qed\ekor

Of course, for an abelian group $G_1$ we have $k=1$ in the above results.
\bem If $w(\bar x,\bar g)=c$ is $2$-largely satisfied in $G^n$, then it is identically satisfied in the abelian quotient $G/G'$. If moreover $G$ is a $BFC$-group, then $G'$ is finite by B.H. Neumann's Lemma \cite{Ne54}, and $G^n$ satisfies a finite disjunction $\bigvee_{c'\in G'}w(\bar x,\bar g)=cc'$.\ebem

We can also deduce results for central elements just from $2$-largeness (although for infinite index $|G_1:Z(G_1)|$ there is no reason that if $X$ is large in $G$ the intersection $X\cap Z(G_1)^n$ is still large in $G\cap Z(G_1)^n$).
\satz If $w(\bar x,\bar g)=c$ is $2$-largely satisfied in $G$, then $w(\bar x,\bar 1)=1$ identically on $G\cap Z(G_1)^n$.\esatz
\bew  Consider $\bar h\in G\cap Z(G_1)^n$. Put $X=\{\bar h'\in G:w(\bar h',\bar g)=1\}$. Then $X\cap \bar h^{-1}X$ is non-empty, so there is $\bar x\in G$ with
$$w(\bar x,\bar g)=c=w(\bar h\cdot\bar x,\bar g)=w(\bar h,\bar 1)\,w(\bar x,\bar g).$$
Hence $w(\bar h,\bar 1)=1$.\ebew
\kor\label{c:10} If $x_1^{k_1}\cdots x_n^{k_n}=c$ is $2$-largely satisfied in $G^n$ and $k=\mbox{\rm gcd}(k_1,\ldots,k_n)$, then $x^k=1$ identically on $Z(G)$.\ekor
\bew We have  $x_1^{k_1}\cdots x_n^{k_n}=1$ on $Z(G)$. Putting $x_i=g\in Z(G)$ and $x_j=1$ for $j\not=i$ we have $g^{k_i}=1$ for all $1\le i\le n$. The result follows.\ebew
\kor If the exponent of $Z(G)$ does not divide $\mbox{\rm gcd}(k_1,\ldots,k_n)$, then $\mu_*(x_1^{k_1}\cdots x_n^{k_n}=c)\le\frac12$.\qed\ekor

\section{Burnside and Engel Equations}\label{s:Burnside}
In Remark \ref{b:exp} we have already seen that if every centraliser of a single element has index at most $k$ in $G$, then $\mu_*(x^m=c)\le1-\frac1{2k}$ unless the exponent of $G$ divides $m$. In this case necessarily $c=x^m=1$.

We shall first prove Miller's Theorem mentioned in the introduction.
\satz Let $c\in G$. If $x^2=c$ is $4$-largely satisfied in $G$, then $G$ is abelian of exponent $2$, and $c=1$.
\esatz
\bew Fix $g,h\in G$. Then there is $x$ with $c=x^2=(gx)^2=(hx)^2=(ghx)^2$. But this implies $x^{-1}gx=g^{-1}$, $x^{-1}hx=h^{-1}$ and $x^{-1}ghx=(gh)^{-1}$. On the other hand,
$$x^{-1}ghx=x^{-1}gx\,x^{-1}hx=g^{-1}h^{-1}=(hg)^{-1}.$$
Hence $gh=hg$ and $G$ is abelian. But now $c=x^2=(gx)^2=g^2x^2=g^2c$, whence $g^2=1$.\ebew
If $G$ satisfies $4$-largely $xax=b$ for some $a,b\in G$, then it satisfies $4$-largely $(ax)^2=ab$, whence $x^2=ab$. Hence $G$ is abelian of exponent $2$, and $a=b$.
\kor If $G$ is not of exponent $2$ or $a\not=b$, then $\mu_*(xax=b)\le\frac34$.\qed\ekor

Recall that the {\em $n^{th}$ Engel condition} is the condition $[x,_n y]=1$, where $[x,_1y]=[x,y]$ and $[x,_{n+1}y]=[[x,_ny],y]$.
Note that 
$$[x,y,y]=[y^{-x}y,y]=y^{-1}y^xy^{-1}y^{-x}yy=[y^{-x},y]^y.$$
Thus the $2$-Engel condition $[x,y,y]=1$ is equivalent to $[y^{-x},y]=1$, that is all conjugacy classes being commutative.
\satz\label{p:18} If $G$ satisfies $7$-largely $x^3=1$ then $G$ is $2$-Engel.\esatz
\bew Put $X=\{g\in G:g^3=1\}$. For $g,h\in G$ consider
$$x\in X\cap g^{-1}X\cap h^{-1}X\cap gX\cap (gh)^{-1}X\cap gh^{-1}X\cap gh^{-1}g^{-1}X.$$
Then $(yx)^3=1$ for $y\in\{1,g,h,g^{-1},gh,hg^{-1},ghg^{-1}\}$, which means that $xyx=y^{-1}x^{-1}y^{-1}$. We calculate the product $xhx^2gx$ in two ways:
$$\begin{aligned}xhx^2gx&=(xhx)(xgx)=h^{-1}(x^{-1}h^{-1}g^{-1}x^{-1})g^{-1}\\
&=h^{-1}ghxghg^{-1}\quad\mbox{and}\\
xhx^2gx&=xh(g^{-1}x)^{-1}x=xh(g^{-1}x)^2x=(xhg^{-1}x)g^{-1}x^2\\
&=gh^{-1}(x^{-1}gh^{-1}g^{-1}x^{-1})=gh^{-1}ghg^{-1}xghg^{-1}.
\end{aligned}$$
Thus $h^{-1}gh=gh^{-1}ghg^{-1}$ and $g^hg=gg^h$. As $h\in G$ was arbitrary, the conjugacy class of $g$ is commutative; as $g$ was arbitrary, all conjugacy classes are commutative.\ebew
\satz Let $G$ be $2$-Engel. If $G$ satisfies $2$-largely $x^3=1$ then $G$ has exponent $3$.\esatz
\bew For any $g\in G$ there is $x\in G$ with $x^3=(gx)^3=1$. As $x^G$ is commutative,
$$g^xg^{-1}g^{x^{-1}}=x^{-1}gxg^{-1}xgx^{-1}=gx^{-g}xx^gx^{-1}=gx^{-g}x^gxx^{-1}=g.$$
Since $g^G$ is commutative, we have
$$\qquad g^3=g^2g^xg^{-1}g^{x^{-1}}=g^2g^{-1}g^{x^{-1}}g^x=(gx)^3=1. \qedhere$$\ebew

\kor\label{c:19} If $G$ satisfies $7$-largely $x^3=1$, then $G$ has exponent $3$. If $G$ is not of exponent $3$ then $\mu_*(x^3=1)\le\frac67$. If moreover $G$ is $2$-Engel, then $\mu_*(x^3=1)\le\frac12$.\qed\ekor

Note that the bound $\frac67$ is not as good as Laffey's bound $\frac 79$ cited in the introduction.
\frag A group which satisfies $5$-largely $x^3=1$, is it $2$-Engel? This would improve our bound to $\frac 45$.\efrag

\kor If $|G:Z(G)|\le 7$ and $G$ satisfies $7$-largely $x^3=c$ for some $c\in G$, then $c=1$ and $G$ has exponent $3$.\ekor
\bew $\{x\in G:x^3=c\}\cap Z(G)$ is $1$-large, whence non-empty, and contains an element $z$. But now there is $x\in G$ with $x^3=1=(zx)^3=z^3x^3=cx^3$, whence $c=1$. We finish by Corollary \ref{c:19}.\ebew

If $|G:Z(G)|$ is prime, then $G$ is abelian, and $2$-largeness is sufficient by Corollary \ref{c:10}.

\section{Commutator Equations}\label{s:[]}
Consider the equation $[x,g]=c$ for some $c,g\in G$. Since $\{x\in G:[x,g]=c\}$ is a coset of $C_G(g)$ or empty, and a coset of a proper subgroup cannot be $2$-large, it follows that if $G$ satisfies $2$-largely $[x,g]=c$ then $g\in Z(G)$ and $c=1$. The following argument generalises this result.
\satz\label{t:13} Suppose $f:G\to H$ satisfies $f(xx')=f(x)^h\,f(x')$ for some $h\in H$ which depends on $x,x'\in G$. If $G_0$ and $G_1$ are groups, $f_0:G_0\to H$ and $f_1:G_1\to H$ are functions such that $G_0\times G\times G_1$ satisfies $k$-largely $f_0(x_0)\,f(x)\,f_1(x_1)=c$ for some $k\ge2$, then $f(G)=1$ and $G_0\times G_1$ satisfies $k$-largely $f_0(x_0)\,f_1(x_1)=c$.\esatz
\bew Fix $g\in G$. By $2$-largeness there is $(x_0,x,x_1)\in G_0\times G\times G_1$ such that
$$f_0(x_0)\,f(x)\,f(x_1)=c=f_0(x_0)\,f(gx)\,f(x_1).$$
Thus $f(x)=f(gx)=f(g)^h\,f(x)$ and $f(g)=1$. It follows that $f_0(x_0)\,f(x)\,f_1(x_1)=f_0(x_0)\,f_1(x_1)$ on $G_0\times G\times G_1$. The result follows.\ebew
\kor\label{c:prod} If $G$ satisfies $2$-largely $\prod_{i<n}[x_i,g_i]=c$ for some $g_i\in G$, then $g_i\in Z(G)$ for all $i<n$ and $c=1$. If not all $g_i$ are central or $c\not=1$ then $\mu_*(\prod_{i<n}[x_i,g_i]=c)\le\frac12$.\ekor
\bew We have $[xx',y]=[x,y]^{x'}[x',y]$.  Now use Theorem \ref{t:13}.\ebew
\bem Theorem \ref{t:13} also holds if $f(xx')=f(x')f(x)^h$, with almost the same proof. Hence Corollary \ref{c:prod} also holds if some factors are of the form $[g_i,x_i]$.\ebem
Gustafson \cite{Gu73} has shown that $\mu_2([x,y]=1)\le\frac12(1+\mu(Z(G))\le\frac58$ for a non-abelian compact topological group $G$, where $\mu$ is the Haar measure on $G$ and $\mu_2$ the product measure on $G^2$. Pournaki and Sobhani \cite{PS08} have generalised this to calculate that $\mu([x,y]=g)<\frac12$ for any $g\not=1$ in a finite group, using Rusin's classification \cite{Ru79} of all finite groups with $\mu([x,y]=1)>\frac{11}{32}$ (see also \cite{FA15}). We have only been able to establish results using $4$-largeness, giving the bound of $\frac34$ in Corollary \ref{k:abel}, so the following two problems remain open:
\frag\begin{enumerate}\item If $G$ satisfies $2$-largely $[x,y]=1$, is $G'=C_2$ and $G/Z(G)$ of exponent $2$, or $G'=C_3$ and $G/Z(G)=S_3$?
\item  If $G$ satisfies $2$-largely $[x,y]=c$ for some $c\in G$, is $c=1$?\end{enumerate}\efrag

\satz If $w(\bar x,\bar g)[x,y]=c$ is satisfied $4$-largely in $G^{n+1}$, where $x\in\bar x$ and $y\notin\bar x$, then $G$ is abelian and $w(\bar x,\bar g)=c$.\esatz
\bew For any $h\in G$ the set 
$$\{(\bar x,x,y):w(\bar x,\bar g)[x,y]=c=w(\bar x,\bar g)[x,hy]\}$$
 is $2$-large in $G^{n+1}$. Hence $\{(x,y)\in G^2:[x,y]=[x,hy]\}$ is $2$-large in $G^2$. Now $[x,hy]=[x,y][x,h]^y$, so $[x,h]=1$ is satisfied $2$-largely in $G$, whence $h\in Z(G)$. It follows that $G$ is abelian. But then $w(\bar x,\bar g)=c$ is satisfied $4$-largely in $G^n$, and must be an identity in $G$ by commutativity and Corollary \ref{k:zentrum}.\ebew
\kor  If $G$ is a group with $\mu_*(w(\bar x,\bar g)[x,y]=c)>\frac34$, then $G$ is abelian satisfying $w(\bar x,\bar g)=c$.\qed\ekor
\kor\label{k:abel} If $G$ satisfies $4$-largely $[x,y]=c$, then $G$ is abelian and $c=1$. If $G$ is not abelian or $c\not=1$, then $\mu_*([x,y]=c)\le\frac34$.\qed\ekor
\bem The same holds for the equation $xcy=yc'x$ with $c\not=c'$: putting $x'=xc$ and $y'=yc'$, this is equivalent to $[x',y']=c^{-1}c'$.\ebem

\satz\label{p:[conj]} Let $g,h\in G$ and $k=\min\{|G:C_G(g)|,|G:C_G(h)|\}$. If $G$ satisfies $k$-largely $[g,h^x]=1$, then $g^G$ and $h^G$ commute.\esatz
\bew If $k=|G:C_G(h)|$, then $\{x\in G:[g,h^x]=1\}\cap C_G(h)$ is $1$-large, whence non-empty, and $[g,h]=1$. Now note that for any $a\in G$ also $|G:C_G(h^a)|=k$ and $[g,h^{ax}]=1$ is satisfied $k$-largely, whence $[g,h^a]=1$ and $[g,h^G]=1$.

If $k=|G:C_G(g)|$, then $\{x\in G:[g^{x^{-1}},h]=1\}\cap C_G(g)$ is $1$-large (still on the left) and non-empty, whence $[g,h]=1$ and we finish as above.\ebew
\kor If $[g^G,h^G]$ is non-trivial for some $g,h\in G$, then $\mu_*([g,h^x]=1)\le 1-\frac 1k$, where $k=\min\{|G:C_G(g)|,|G:C_G(h)|\}$.\qed\ekor

\satz If $g,h,c\in G$ and $[x,g,h]=c$ is $2k$-largely satisfied, where $k=|G:C_G(h)|$, then $[G,g,h]=1$.
Similarly, if $[g,x,h]=c$ is $2k$-largely satisfied for some $c\in Z(G)$, then $[g,G,h]=1$.\esatz
\bew Choose $a\in G$. Then the set $X=\{x\in G:[x,g,h]=c=[ax,g,h]\}$ is $k$-large, and for $x\in X$ we have
$$[x,g,h]=c=[ax,g,h]=[[a,g]^x[x,g],h]=[[a,g]^x,h]^{[x,g]}[x,g,h],$$
whence $[[a,g]^x,h]=1$. By Theorem \ref{p:[conj]} we have $[a,g,h]=1$.

If $[g,x,h]=c$ is $2k$-largely satisfied with $c\in Z(G)$, then for $a\in G$ we obtain a $k$-large $X\subseteq G$ such that for $x\in X$ we have
$$[g,x,h]=c=[g,ax,h]=[[g,x][g,a]^x,h]=[g,x,h]^{[g,a]^x}[[g,a]^x,h],$$
whence $[[g,a]^x,h]=1$, and $[g,a,h]=1$ by Theorem \ref{p:[conj]}.\ebew
 
\kor If $g,h\in G$ and $k=|G:C_G(h)|$, then $[G,g,h]\not=1$ implies $\mu_*([x,g,h]=c)\le1-\frac1{2k}$ for any $c\in G$, and 
$[g,G,h]\not=c$ implies $\mu_*([g,x,h]=c)\le1-\frac1{2k}$ for any $c\in Z(G)$.\qed\ekor

We shall now generalise Corollary \ref{k:abel} to higher nilpotency classes. However, the proof requires an additional assumption.
\satz\label{t:nilp} Suppose $s<\omega$ is such that for all $i<k$ there is a set $A_i$ of size at most $s$ such that  $Z(G/Z_i(G))=C_{G/Z_i(G)}(A_i)$. If $G$ satisfies $2(s+1)^k$-largely $[x_0,x_1,\ldots,x_k]=c$, then $c=1$ and $G$ is nilpotent of class at most $k$.
\esatz
\bew We use induction on $k$. For $k=1$ note that $s\ge 1$ (otherwise $G$ is abelian and we are done), so the result follows from Corollary \ref{k:abel}.

Now suppose the assertion is true for $k$, and 
$$X=\{\bar x\in G^{k+2}:[x_0,x_1,\ldots,x_{k+1}]=c\}$$ is $2(s+1)^{k+1}$-large in $G^{k+2}$. If $A_0=\{a_i:i<s\}$ consider the projection $Y$ of $X\cap\bigcap_{i<s}(1,\ldots,1,a_i^{-1})X$ to the first $k+1$ coordinates, and note that it is $2(s+1)^k$-large. Then for all $(x_0,\ldots,x_k)\in Y$ there is $y\in G$ such that
$$[x_0,\ldots,x_k,y]=c=[x_0,\ldots,x_k,a_iy]=[x_0,\ldots,x_k,y]\,[x_0,\ldots,x_k,a_i]^y$$
for all $i<s$, whence $[x_0,\ldots,x_k]\in Z(G)$. By inductive assumption $G/Z(G)$ is nilpotent of class at most $k$, and we are done.\ebew

\kor\label{k:nilp} Let $s$ be as above. If $G$ is not nilpotent of class at most $k$ or $c\not=1$, then $\mu_*([x_0,x_1,\ldots,x_k]=c)\le1-\frac12(s+1)^{-k}$.\qed\ekor
\bem Recall that an \Mc-group is a group $G$ such that for every subset $A$ there is a finite subset $A_0\subseteq A$ such that $C_G(A)=C_G(A_0)$. Equivalently, $G$ satisfies the ascending (or the descending) chain condition on centralisers. Roger Bryant \cite{Br79} has shown that in an \Mc-group, for every iterated centre $Z_i(G)$ there is a finite set $A_i$ such that $Z(G/Z_i(G))=C_{G/Z_i(G)}(A_i)$. So in an \Mc-group we can find some $s$ as needed for Theorem \ref{t:nilp} and Corollary \ref{k:nilp}.\ebem
\frag To what extent do we need the \Mc-condition (or similar) in Theorem \ref{t:nilp} and Corollary \ref{t:nilp}? 
It is not needed for nilpotency class $1$ (Corollary \ref{k:abel}). In general, assuming just $2^{k+1}$-largeness of $[x_0,\ldots,x_k]=c$, we obtain that $\{\bar x\in G^k:[x_0,\ldots,x_{k-1}]\in C_G(g)\}$ is $2^k$-large in $G^k$ for any $g\in G$. Does this imply $\gamma_k(G)\le C_G(g)$, or even $\gamma_k(G)\le Z(G)$?\efrag

\section{Nilpotent groups}\label{s:Nilp}
We shall first introduce the notion of a supercommutator from \cite{JW00}.
\defn Any variable and any constant from $G$ is a {\em supercommutator}; if $v$ and $w$ are supercommutators, then $v^{-1}$ and $[v,w]$ are supercommutators.\edefn
Alternatively, we could have said that $x$, $x^{-1}$ and $g$ are supercommutators for any variable $x$ and any $g\in G$, and that if $v$ and $w$ are supercommutators, so is $[v,w]$.
\defn The set $\Var(v)$ of variables of a supercommutator $v$ is defined by $\Var(x)=\{x\}$, $\Var(g)=\emptyset$, $\Var(v^{-1})=\Var(v)$, and $\Var([v,w]=\Var(v)\cup\Var(w)$. We put $\var(v)=|\Var(v)|$, the {\em variable number} of $v$. If $\bar x$ is a tuple of variables, we put $\Var_{\bar x}=\Var(v)\cap\bar x$, $\Var'_{\bar x}(v)=\Var(v)\setminus\bar x$, $\var_{\bar x}(v)=|\Var_{\bar x}(v)|$ and $\var'_{\bar x}(v)=|\Var'_{\bar x}(v)|$.
%The {\em complexity} $\comp(v)$ of $v$ is $1$ if $v$ is a variable and $0$ if $v$ is a constant in $H$; 
%$\comp(v^{-1})=\comp(\sigma(v))=\comp(v)$, and $\comp([v,w])=\comp(v)+\comp(w)$.
\edefn
Clearly $\var([v,v'])\ge\max\{\var(v),\var(v')\}$, and similarly for $\var_{\bar x}$ and $\var'_{\bar x}$.
\lem\label{product} Let $H\trianglelefteq G$ and $v(\bar x,\bar z)$ a supercommutator.\begin{enumerate}
%\item\label{one} $\var_{\bar x}(v)\le \var(v)\le \comp(v)$.
\item\label{two} $v$ defines a function from $H^{|\bar x\bar z|}$ to $\gamma_{\var(v)}(H)$.
\item\label{three} If $\var_{\bar x}(v)>0$ and $\bar x$, $\bar y$ and $\bar z$ are pairwise disjoint, then
$$v(\bar y\cdot\bar x,\bar z)=v(\bar x,\bar z)\,v(\bar y,\bar z)\,\Phi(\bar x,\bar y,\bar z),$$
where $\Phi$ is a product of supercommutators whose factors $w$ satisfy\begin{itemize}
\item[(\dag)] $\Var_{\bar z}(w)=\Var_{\bar z}(v)$, and if $x_i\in\Var_{\bar x}(v)$ then $x_i\in\Var(w)$ or $y_i\in\Var(w)$, and both possibilities occur for at least one~$i$.\end{itemize}
\item\label{four} If $v(\bar x,\bar z)$ is a product of supercommutators whose factors $w$ satisfy $\var_{\bar x}(w)>0$ and $\var'_{\bar x}(w)\ge n$, then
$$v(\bar y\cdot\bar x,\bar z)=v(\bar x,\bar z)\,v(\bar y,\bar z)\,\Phi(\bar x,\bar y,\bar z),$$
where $\Phi$ is a product of supercommutators whose factors $w$ satisfy  $\var_{\bar x}(w)>0$ and $\var'_{\bar x}(w)>n$.
\end{enumerate}\elem
\bew $(\ref{two})$ is proved as in \cite[Lemme 6(1)]{JW00} by induction, using that $\gamma_n(H)$ is characteristic in $H$, whence normal in $G$, and $[\gamma_n(H),\gamma_m(H)]\le\gamma_{n+m}(H)$. We shall show $(\ref{three})$ by induction on the construction of $v$.

If $v=x\in\bar x$ we have $v(yx)=yx=xy[y,x]=v(x)v(y)[y,x]$; if $v=x^{-1}$ we have $v(yx)=x^{-1}y^{-1}=v(x)v(y)$. This leaves the case $v=[v_1,v_2]$ for two supercommutators $v_1$ and $v_2$. We shall assume $\var_{\bar x}(v_1)>0$ and $\var_{\bar x}(v_2)>0$ (the case $\var_{\bar x}(v_1)\var_{\bar x}(v_2)=0$ is analogous, but simpler).

By inductive hypothesis, there are $\Phi_i$ for $i=1,2$, products of supercommutators satisfying (\dag) relative to $v_i$, such that $$v_i(\bar y\cdot\bar x,\bar z)=v_i(\bar x,\bar z)\,v_i(\bar y,\bar z)\,\Phi_i.$$
Then
$$\begin{aligned}v(\bar y\cdot\bar x,\bar z)&=[v_1(\bar y\cdot\bar x,\bar z),v_2(\bar y\cdot\bar x,\bar z)]\\
&=[v_1(\bar x,\bar z)\,v_1(\bar y,\bar z)\,\Phi_1,v_2(\bar x,\bar z)\,v_2(\bar y,\bar z)\,\Phi_2]\\
&=[v_1(\bar x,\bar z),v_2(\bar x,\bar z)]\,[v_1(\bar y,\bar z),v_2(\bar y,\bar z)]\,\Phi
=v(\bar x,\bar z)\,v(\bar y,\bar z)\,\Phi,\end{aligned}$$
where $\Phi$ is a product of supercommutators $[w,w']$\begin{enumerate}
\item[(i)] where $w\in\Phi_1\cup\{v_1(\bar x,\bar z),v_1(\bar y,\bar z)\}$ and
$w'\in\Phi_2\cup\{v_2(\bar x,\bar z),v_2(\bar y,\bar z)\}$, except for $[v_1(\bar x,\bar z),v_2(\bar x,\bar z)]$ and $[v_1(\bar y,\bar z),v_2(\bar y,\bar z)]$; it is clear that these must satisfy (\dag).
\item[(ii)] where one of $w,w'$ is from (i), so $[w,w']$ satisfies (\dag).
\item[(iii)] where one of $w,w'$ is equal to $v(\bar x,\bar z)$ and the other contains at least one $y_i$, or one is equal to $v(\bar y,\bar z)$ and the other contains at least one $x_i$; again $[w,w']$ satisfies (\dag).
\item[(iv)] which are obtained iteratively from supercommutators from (ii) and (iii) by commutation with other supercommutators, thus satisfying (\dag).\end{enumerate}
Here (i) takes care of the commutators of various factors of the two products, while (ii)--(iv) takes care of the correct order. Note that the only factor without a variable $y_i$ is $v(\bar x,\bar z)$ and the only factor without a variable $x_j$ is $v(\bar y,\bar z)$. 

To show $(\ref{four})$ note first that for a single supercommutator $v$ the factorisation given in $(\ref{three})$ satisfies the requirement. So for a product of supercommutators, we apply $(\ref{three})$ to every factor, and then use commutators to get them into the right order. Note that we never have to commute a $w(\bar x,\bar z)$ with a $w'(\bar x,\bar z)$, or a $w(\bar y,\bar z)$ with a $w'(\bar y,\bar z)$, as they already appear in the correct order with respect to one another. It follows that all new commutators satisfy (\dag), whence $\var'_{\bar x}>n$.\ebew
\satz\label{constant} If $G$ is nilpotent of class $k$ and $v$ is a product of supercommutators $w$ with $\var_{\bar x}(w)>0$ and $\var'_{\bar x}(w)\ge n$ such that $G$ satisfies $\max\{2^{k-n},1\}$-largely $v(\bar x,\bar g)=c$, then $c=1$.\esatz
\bew This  is true for $n\ge k$, as then $\var(w)=\var_{\bar x}(w)+\var'_{\bar x}(w)\ge 1+n$, and
$$c=w(\bar x,\bar g)\in\gamma_{\var(w)}G\le\gamma_{n+1}G=\{1\}$$
for some $\bar x\in G$.

Now suppose it is true for $n+1\le k$, and let $v(\bar x,\bar z)$ be a product of supercommutators $w$ with $\var_{\bar x}(w)>0$ and $\var'_{\bar x}\ge n$, such that $H$ satisfies $2^{k-n}$-largely $v(\bar x,\bar g)=c$. By Lemma \ref{product} there is $\Phi$, a product of supercommutators whose factors $w$ satisfy $\var_{\bar x}(w)>0$ and $\var'_{\bar x}(w)>n$, such that
$$v(\bar y\cdot\bar x,\bar z)=v(\bar x,\bar z)\,v(\bar y,\bar z)\,\Phi(\bar x,\bar y,\bar z).$$
Choose $\bar h\in G$ with $v(\bar h,\bar g)=c$. If $X=\{\bar x\in G:v(\bar x,\bar g)=c\}$, then $X$ is $2^{k-n}$-large, and $Y=X\cap\bar h^{-1}X$ is $2^{k-n-1}$-large. Moreover, for $\bar x\in Y$ we have
$$\Phi(\bar x,\bar h,\bar g)=v(\bar h,\bar g)^{-1}v(\bar x,\bar g)^{-1}v(\bar h\cdot\bar x,\bar g)=c^{-1}c^{-1}c=c^{-1}.$$
By hypothesis $c^{-1}=1$ and we are done.\ebew
\satz\label{nilpotent} If $G$ is nilpotent of class $k$ and satisfies $2^k$-largely an equation $v(\bar x,\bar g)=c$, then it satisfies $v(\bar x,\bar g)=c$.\esatz
\bew Bringing all the constants to the right-hand side, we may assume that $v(\bar x,\bar z)$ is a product of supercommutators $w$ with $\var_{\bar x}(w)>0$. By Lemma \ref{product} there is $\Phi$, a product of supercommutators whose factors $w$ satisfy $\var_{\bar x}(w)>0$ and $\var'_{\bar x}(w)>0$, such that
$$v(\bar y\cdot\bar x,\bar z)=v(\bar x,\bar z)\,v(\bar y,\bar z)\,\Phi(\bar x,\bar y,\bar z).$$
Fix $\bar h\in G$. Then
$$\Phi(\bar x,\bar h,\bar g)=v(\bar h,\bar g)^{-1}c^{-1}c=v(\bar h,\bar g)^{-1}$$
$2^{k-1}$-largely on $G$. By Theorem \ref{constant} we have $v(\bar h,\bar g)=1$. So $v(\bar x,\bar g)$ is constant.\ebew
\kor If $G$ is nilpotent of class $k$ and $x^n=c$ is true $2^k$-largely, then $c=1$ and the exponent of $G$ divides $n$.\ekor
\bew Immediate from Theorem \ref{nilpotent}.\ebew
\kor If $G$ is nilpotent of class $k$ and $\mu_*(x^n=c)>1-2^{-k}$, then $c=1$ and the exponent of $G$ divides $n$.\qed\ekor

\section{Autocommutativity}\label{s:Auto}
The notion of autocommutativity has been introduced by Sherman in 1975 \cite{Sh75}.
\defn Let $G$ be a finite group, $\Sigma$ a group of automorphisms of $G$, and $H$ a subgroup of $G$. The {\em degree of autocommutativity relative to $(H;\Sigma)$} is given by
$$\ac(H;\Sigma)=\frac{|\{(\sigma,g)\in\Sigma\times H:\sigma(g)=g\}|}{|\Sigma|\cdot|H|}.$$
It gives the probability that a random element of $H$ is fixed by a random automorphism in $\Sigma$.\edefn
Note that $\ac(H;\Sigma)=\mu(\{(\sigma,g)\in\Sigma\times H:\sigma(g)=g\})$, where $\mu$ is the counting measure on $\Sigma\times H$.
\satz\label{t:ac} Let $H\le G$ be finite groups, $\Sigma$ a group of automorphisms of $G$, and suppose that $\{(\sigma,g)\in\Sigma\times H:\sigma(g)=g\}$ is $4$-large in $\Sigma\times H$. Then $H\le\mathrm{Fix}(\Sigma)$.\esatz
\bew Given $\sigma\in\Sigma$ and $g\in H$, by $4$-largeness there are $x\in H$ and $\tau\in\Sigma$ with
$$\tau(x)=x,\quad(\sigma\circ\tau)(x)=x,\quad\tau(gx)=gx\quad\mbox{and}\quad(\sigma\circ\tau)(gx)=gx.$$
Then
$$gx=\sigma(\tau(gx))=\sigma(gx)=\sigma(g)\sigma(x)=\sigma(g)\sigma(\tau(x))=\sigma(g)x,$$
whence $g=\sigma(g)$.\ebew
\kor If $H\le G$ are finite groups and $\Sigma$ is a group of automorphisms of $G$ with $H\not\le\mbox{Fix}(\Sigma)$, then
$\ac(H;\Sigma)\le\frac34$.\ekor
\bew If $\ac(H;\Sigma)>\frac34$ then $\{(\sigma,g)\in\Sigma\times H:\sigma(g)=g\}$ is $4$-large in $\Sigma\times H$ by Lemma \ref{l:mu}. Hence $H\le\mathrm{Fix}(\Sigma)$ by Theorem \ref{t:ac}.
\ebew

\end{document}